 \documentclass[twoside,reqno,11pt]{fcaa-var}
 \usepackage{amsfonts}
 \usepackage{amsfonts}
 \usepackage{amsfonts}
 \usepackage{cases}
 \usepackage{amsmath}
 \usepackage{amsmath}

\usepackage{graphicx}
\usepackage{epsfig}
\usepackage{amsthm}
\usepackage{amsmath}
\usepackage{latexsym}
\usepackage{amsfonts}
\usepackage{amssymb}

\newtheoremstyle{theorem}
  {15pt}          
  {15pt}  
  {\sl}  
  {\parindent}
  {\sc}  
  {. }   
  { }    
  {}     
\theoremstyle{theorem}

\newtheoremstyle{defi}
  {15pt}          
  {15pt}  
  {\rm}  
  {\parindent}     
  {\sc}  
  {. }    
  { }    
  {}     
\theoremstyle{defi}

\newtheorem{remark}{Remark}[section]

 \def\proofname{\indent\rm P\ r\ o\ o\ f.}
 
 \def\qed{\hfill$\Box$} 

 \usepackage{hyperref}
 \def\theequation{\arabic{section}.\arabic{equation}}

\usepackage{threeparttable}
\usepackage[figuresright]{rotating}
 \usepackage{multirow}

  \def\themycopyrightfootnote{
       \hspace*{-0.7cm}
    {\normalsize  \copyright \, 2015\,  Diogenes  Co., Sofia  \\
        \noindent  pp. 1492--1506\,, DOI: 10.1515/fca-2015-0086
     \hfill  \vspace*{-0.45cm}
   \hspace*{10.4cm}  
    }}

 \title[Global Pad\'{e} approximations of \dots]
 {Global Pad\'{e} approximations of \\ [3pt]
  the generalized Mittag-Leffler function \\ [3pt] and its inverse}
 \author[\normalsize C. Zeng, Y.Q. Chen]{\normalsize Caibin Zeng $^1$, YangQuan Chen $^2$}

\begin{document}


\bigskip

 \noindent \textbf{RESEARCH PAPER}

 \bigskip

 \begin{abstract}

This paper proposes a global Pad\'{e} approximation of
the generalized Mittag-Leffler function $E_{\alpha,\beta}(-x)$ with
$x\in[0,+\infty)$. This uniform approximation can account for both
the Taylor series for small arguments and asymptotic series for
large arguments. Based on the complete monotonicity of the function
$E_{\alpha,\beta}(-x)$, we work out the global Pad\'{e}
approximation [1/2] for the particular cases $\{0<\alpha<1,
\beta>\alpha\}$, $\{0<\alpha=\beta<1\}$, and $\{\alpha=1,
\beta>1\}$, respectively. Moreover, these approximations are
inverted to yield a global Pad\'{e} approximation of the inverse
generalized Mittag-Leffler function $-L_{\alpha,\beta}(x)$ with
$x\in(0,1/\Gamma(\beta)]$. We also provide several examples with
selected values $\alpha$ and $\beta$ to compute the relative error from the
approximations. Finally, we point out the possible applications
using our established approximations in the ordinary and partial
time-fractional differential equations in the sense of
Riemann-Liouville.

 \medskip


 \smallskip


 \end{abstract}

 \maketitle

 \vspace*{-16pt}



\section{Introduction}
\label{sec:1}

\setcounter{section}{1}
\setcounter{equation}{0}\setcounter{theorem}{0}

The Mittag-Leffler function and its generalizations are very
important special functions that find widespread use in the framework of
fractional calculus. Similarly to the exponential function frequently
used in the solutions of integer-order differential equations, the (generalized)
Mittag-Leffler functions play an analogous role in the solution of
fractional-order differential equations. In fact, the exponential function itself
is of a very specific form, one of an infinite set, of these seemingly
ubiquitous functions. The standard definition of Mittag-Leffler
function was given by \cite{Mittag1903}
\begin{equation}\label{eq1}
    E_\alpha(z)=\sum_{k=0}^{\infty}\frac{z^k}{\Gamma(\alpha k+1)},
\, ~~ z\in\mathbb{C},
\end{equation}
where $\Re (\alpha)>0$. The Mittag-Leffler function with two parameters
(sometimes also called the generalized Mittag-Leffler function),
appeared most frequently and had the following form \cite{Wiman1905}
\begin{equation}\label{eq2}
    E_{\alpha,\beta}(z)=\sum_{k=0}^{\infty}\frac{z^k}{\Gamma(\alpha k+\beta)},
\, ~~ z\in\mathbb{C},
\end{equation}
where $\Re (\alpha)>0$ and $\beta \in \mathbb{C}$. For $\beta=1$, we have
$E_\alpha(z)=E_{\alpha,1}(z)$, and also, $E_{1,1}(z)=e^z$. Later on,
some studies
\cite{Erdelyi1955,Dzherbashyan1966,Kilbas1996,Mainardi2000,Kilbas2004,Peng2010,Haubold2011}
contributed to several properties and applications of the
(generalized) Mittag-Leffler functions. Nowadays the Mittag-Leffler
function is referred to as the {\it Queen Function of Fractional Calculus}
\cite{Mainardi2007}. In this respect we also recommend some
classical books on fractional calculus
\cite{Podlubny1999,Kilbas2006,Mainardi2010,Sheng2012}.

There are several approaches to the numerical computation
of the (generalized) Mittag-Leffler function.
Gorenflo et al. \cite{Gorenflo2002}
studied the computation of the generalized Mittag-Leffler function
$E_{\alpha,\beta}(z)$ and its derivative for all values of
the parameters $\alpha>0$, $\beta\in\mathbb{R}$, and all values of the
argument $z\in\mathbb{C}$. Hilfer and Seybold \cite{Hilfer2006,Seybold2008} reported a
simpler algorithm to compute the generalized Mittag-Leffler function
based on mixed techniques such as Taylor series,
asymptotic series, and integral representations.
It is also worth pointing out
that Podlubny \cite{Podlubny2005} provided a MATLAB routine for
evaluating the generalized Mittag-Leffler function with desired
accuracy. Moreover, Garrappa \cite{Garrappa2015a} presented an efficient method to
evaluate the Mittag-Leffler function
based on the numerical inversion of its Laplace transform.
The corresponding MATLAB code is made freely available, \cite{Garrappa2015b}.

Recently, Starovotov and Starovotova \cite{Starovoitov2007}
discussed the Pad\'{e} approximations for the Mittag-Leffler
function and showed that the approximations serve uniformly on the
compact set $\{|z| \le 1\}$. Mainardi \cite{Mainardi2013} also used
the Pad\'{e} approximation to provide lower and upper bounds to the
Mittag-Leffler function $E_{\alpha}(-t^\alpha)$ for $t > 0$. The
Pad\'{e} approximation is better than a truncated Taylor series, which is not necessarily compatible with the asymptotic behavior for
large arguments. Concezzi and Spigler \cite{Concezzi2015} proved part of Mainardi's
conjecture by simple classical estimates.
Diethelm et al. \cite{Diethelm2005} provided a table of coefficients
of the rational approximants to the function $E_{\alpha}(-t^\alpha)$
for $0<\alpha<1$. The coefficients of these approximations were
numerically calculated to approximate the Mittag-Leffler
function on the interval $[0.1, 15]$ with a certain accuracy.
Indeed, Winitzki \cite{Winitzki2003} provided
the so-called global Pad\'{e} approximation, which constructs
uniform approximations to analytic transcendental functions. These
uniform approximations are built from elementary functions using
both Taylor and asymptotic series of the given transcendental
function. The authors applied it to find good approximations of
several functions including the elliptic function, the error
function of real and imaginary arguments, the Bessel functions, and
the Airy function. Also, Atkinson and Osseiran \cite{Atkinson2011}
applied this method to find a uniform rational approximation of the
Mittag-Leffler function $E_\alpha(-z)$ with $0<\alpha<1$ and
$z\in(0,\infty)$.

However, to the best of our knowledge, a
global Pad\'{e} approximation is need for the generalized Mittag-Leffler
function $E_{\alpha,\beta}(z)$ for all parameters $\alpha$, $\beta$
and $z$.

On the other hand, the completely monotone functions are known to
play an important role in different branches of mathematics and
especially in the probability theory (see \cite{Feller1971} for
example). Recall that a function $f(x)$ is said to be completely
monotonic on an interval $I$ if $f(x)$ has
derivatives of all orders on $I$ which alternate successively
in sign, that is
\begin{equation}\label{eq0}
    (-1)^m f^{(m)}\ge0
\end{equation}
for $x \in I$ and $m\ge 0$.
In particular, Pollard \cite{Pollard1948} proved that the
Mittag-Leffler function $E_\alpha(-x)$ with $x\ge0$ is completely
monotonic, that is
\begin{equation*}
    (-1)^m \frac{d^m}{dx^m}E_\alpha(-x)\ge0
\end{equation*}
for all $m=0, 1,2,\ldots$, if $0\le \alpha \le 1$. Based on the corresponding probability measures and the Hankel
contour integration, Schneider \cite{Schneider1996} proved that the
generalized Mittag-Leffler function $E_{\alpha,\beta}(-x)$ with
$x\ge0$ is completely monotonic if and only if $0<\alpha\le 1$ and
$\beta\ge\alpha$. In other word, it yields
\begin{equation}\label{eq3}
    (-1)^m \frac{d^m}{dx^m}E_{\alpha,\beta}(-x)\ge0
\end{equation}
for all $m=0, 1,2,\ldots$, if $0< \alpha \le 1,
\beta\ge\alpha$. This result was also proved in a simpler way
\cite{Miller1999}. This property is essential for the discussion of
the inverse generalized Mittag-Leffler function below.

We will focus on the generalized Mittag-Leffler function
$E_{\alpha,\beta}(-x)$ and its inverse with $0<\alpha\le 1$ and
$\beta\ge\alpha$. More precisely, we will divide the paper into three
particular cases $\{0<\alpha<1, \beta>\alpha\}$,
$\{0<\alpha=\beta<1\}$, and $\{\alpha=1, \beta>1\}$, for which one
can find good approximations by using both Taylor and asymptotic
series.

The paper is organized as follows. In Section \ref{sec:2} we develop the
global Pad\'{e} approximations of the generalized Mittag-Leffler
function $E_{\alpha,\beta}(-x)$ for particular cases $\{0<\alpha<1,
\beta>\alpha\}$, $\{0<\alpha=\beta<1\}$, and $\{\alpha=1,
\beta>1\}$, respectively. In Section \ref{sec:3} we study the inverse
generalized Mittag-Leffler function $-L_{\alpha,\beta}(x)$ and
obtain its uniform approximations. Finally, we give some concluding
discussions in Section \ref{sec:4} and close the paper.


\section{Global Pad\'{e} approximations of $E_{\alpha,\beta}{(-x)}$}
\label{sec:2}

\setcounter{section}{2}
\setcounter{equation}{0}\setcounter{theorem}{0}

We now consider the global Pad\'{e} approximations of the
generalized Mittag-Leffler function $E_{\alpha,\beta}(-x)$ in the
domain $[0, +\infty)$. Note that the function $E_{\alpha,\beta}(-x)$
is finite everywhere in the interval $[0, +\infty)$ for $0<\alpha\le
1$ and $\beta\ge\alpha$. So we can apply the idea for nonsingular
functions by Winitzki \cite{Winitzki2003} to obtain a global
Pad\'{e} approximation of the generalized Mittag-Leffler function.
Also, we notice the fact that $E_{\alpha,\beta}(+\infty)=0$, which
motivates us to choose the series at $x=+\infty$ starting with a
higher power of $x^{-1}$. Additionally, we adopt the relative error
to measure the discrepancy between the exact value and the
global Pad\'{e} approximation. Given the exact value $v_{exact}(x)\neq0$ of the function
$E_{\alpha,\beta}(-x)$ and its global Pad\'{e} approximation
$v_{approx}(x)$, the relative error is
$\eta=(\eta_1,\ldots,\eta_k,\ldots,\eta_n)$, where
\begin{equation}\label{eq4}
    \eta_k=\left|\frac{v_{exact}(x_k)-v_{approx}(x_k)}{v_{exact}(x_k)}\right|\times
    100\%,
\end{equation}
with the partition $P=\{0\le x_1<\cdots<x_k<\cdots<x_n<\infty\}$ of variable $x$.

For the case $\{0<\alpha<1, \beta>\alpha\}$, it follows from the
definition (\ref{eq2}) that
\begin{equation}\label{eq5}
\left\{
\begin{aligned}
    \Gamma(\beta-\alpha)x E_{\alpha,\beta}(-x)&=
    \Gamma(\beta-\alpha)x \sum_{k=0}^{m-2}\frac{(-x)^k}
    {\Gamma(\beta+\alpha k)}+\mathcal
    {O}(x^m)\\
    &\equiv a(x)+ \mathcal
    {O}(x^m),\\
    \Gamma(\beta-\alpha)x E_{\alpha,\beta}(-x)&=
    -\Gamma(\beta-\alpha)x \sum_{k=1}^{n}\frac{(-x)^{-k}}{\Gamma(\beta-\alpha k)}+\mathcal
    {O}(x^{-n})\\
    &\equiv b(x^{-1})+ \mathcal
    {O}(x^{-n}),
\end{aligned}
\right.
\end{equation}
at $x = 0$ and $x = +\infty$, respectively. The multiplication by
$\Gamma(\beta-\alpha)x$ ensures that the first coefficient of the
asymptotic series (\ref{eq5}) is 1.

We now look for a rational approximation of the form
\begin{equation}\label{eq6}
    \Gamma(\beta-\alpha)x E_{\alpha,\beta}(-x)\approx \frac{p(x)}{q(x)}=\frac{p_0+p_1x+\ldots+p_\nu x^\nu}{q_0+q_1x+\ldots+q_\nu
    x^\nu},
\end{equation}
where $\nu$ is an appropriately chosen integer. The problem is to
find the coefficients $p_i$ and $q_i$ such that (\ref{eq6}) has the
correct expansions at $x = 0$ and $x = +\infty$. Since the leading
term of (\ref{eq6}) at $x = +\infty$ is $p_\nu/q_\nu$, we can set
$p_\nu=q_\nu = 1$. This formulation is similar to the problem of
Hermite-Pad\'{e} interpolation with two anchor points, except that
one of the points is at infinity where we use an expansion in
$x^{-1}$. The unknown coefficients $p_i$ and $q_i$ are found from
the system of linear equations written compactly as
\begin{equation}\label{eq7}
    p(x)-q(x)a(x)=\mathcal
    {O}(x^m)~\text{at}~x=0,
\end{equation}
\vspace*{-8pt}
\begin{equation}\label{eq8}
    \frac{p(x)}{x^\nu}- \frac{q(x)}{x^\nu}b(x^{-1})=\mathcal
    {O}(x^{-n})~\text{at}~x=+\infty.
\end{equation}
Here it is implied that the surviving polynomial coefficients in $x$
or $x^{-1}$ are same for the both sides of (\ref{eq5}). This assumes that $p(x)$ and $q(x)$ have no
common polynomial factors. Moreover, these two equations form an
inhomogeneous linear system of $(m+n-1)$ equations for $2\nu$
unknowns $p_i$, $q_i$, $0\le i \le \nu-1$. So, when a solution
exists, it is unique if $m+n$ is odd.

By setting $\nu = 2$, we thus search for an approximation for the
function $\Gamma(\beta-\alpha)x E_{\alpha,\beta}(-x)$ of the form
\vskip -12pt
\begin{equation}\label{eq9}
    \Gamma(\beta-\alpha)x E_{\alpha,\beta}(-x)\approx \frac{p_0+p_1x+x^2}{q_0+q_1x+x^2},
\end{equation}
\vskip -3pt \noindent
where the two series in (\ref{eq5}) are truncated to
orders $m = 3$ and $n = 2$, yielding the functions
\vskip -13pt
\begin{equation}\label{eq10}
    a(x)=\frac{\Gamma(\beta-\alpha)}{\Gamma(\beta)}x-\frac{\Gamma(\beta-\alpha)}{\Gamma(\beta+\alpha)}x^2,
\end{equation}
\vspace*{-8pt}
\begin{equation}\label{eq11}
    b(x^{-1})=1-\frac{\Gamma(\beta-\alpha)}{\Gamma(\beta-2\alpha)}x^{-1}.
\end{equation}
Substituting (\ref{eq9})-(\ref{eq11}) into (\ref{eq7}) and
(\ref{eq8}) and collecting equal powers of $x$ through $\mathcal
    {O}(x^m)$ and $\mathcal{O}(x^{-n})$, we obtain
    \vskip -10pt
\begin{equation*}
 \left\{
\begin{array}{l}
 p_0=0,\\[2pt]
 p_1-\frac{\Gamma(\beta-\alpha)}{\Gamma(\beta)}q_0=0,\\ [2pt]
 1+\frac{\Gamma(\beta-\alpha)}{\Gamma(\beta+\alpha)}q_0-\frac{\Gamma(\beta-\alpha)}{\Gamma(\beta)}q_1=0,\\ [2pt]
 p_1-q_1+\frac{\Gamma(\beta-\alpha)}{\Gamma(\beta-2\alpha)}=0.
 \end{array} \right.
\end{equation*}
\vskip -3pt \noindent
Its solution can be expressed by
\vskip -10pt
\begin{equation}\label{eq12}
\left\{
\begin{array}{l}
 p_0=0,\\ [2pt]
 p_1=\frac{\Gamma(\beta)\Gamma(\beta+\alpha)-\frac{\Gamma(\beta+\alpha)\Gamma(\beta-\alpha)^2}
    {\Gamma(\beta-2\alpha)}}{\Gamma(\beta+\alpha)\Gamma(\beta-\alpha)-\Gamma(\beta)^2}\,,\\ [3pt]
 q_0=\frac{\frac{\Gamma(\beta)^2\Gamma(\beta+\alpha)}{\Gamma(\beta-\alpha)}-
    \frac{\Gamma(\beta)\Gamma(\beta+\alpha)\Gamma(\beta-\alpha)}
    {\Gamma(\beta-2\alpha)}}{\Gamma(\beta+\alpha)\Gamma(\beta-\alpha)-\Gamma(\beta)^2}\,,\\ [3pt]
 q_1=\frac{\Gamma(\beta)\Gamma(\beta+\alpha)-\frac{\Gamma(\beta)^2\Gamma(\beta-\alpha)}
    {\Gamma(\beta-2\alpha)}}{\Gamma(\beta+\alpha)\Gamma(\beta-\alpha)-\Gamma(\beta)^2}\,.
 \end{array}\right.
\end{equation}
Therefore, we obtain a global Pad\'{e} approximation from
(\ref{eq9}) and (\ref{eq12}), 
\vskip -10pt
\begin{eqnarray}\label{eq13}
   E_{\alpha,\beta}(-x)  \approx \frac{1}{\Gamma(\beta-\alpha)x}\frac{p_0+p_1x+x^2}{q_0+q_1x+x^2}
   =
   \frac{\frac{1}{\Gamma(\beta)}+\frac{1}{\Gamma(\beta-\alpha)q_0}x}{1+\frac{q_1}{q_0}x+\frac{1}{q_0}x^2}\,.
\end{eqnarray}

As an example, for $\alpha = 1/2$ and $\beta=3/2$, it follows that
\begin{equation*}
    E_{\frac{1}{2},\frac{3}{2}}(-x)=\left\{
                                      \begin{array}{ll}
                                        \frac{2}{\sqrt{\pi}}, & \hbox{if $x=0$,} \\[4pt]
                                        \frac{1-\exp(x^2)\text{erfc}(x)}{x}, & \hbox{if $x>0$,}
                                      \end{array}
                                    \right.
\end{equation*}
due to the properties $E_{{1}/{2}}(z)=\exp(z^2)\text{erfc}(-z)$,
~$E_{\alpha,\beta}(z)=zE_{\alpha,\alpha+\beta}(z)+1/\Gamma(\beta)$,
where erfc$(x)$ is complementary error function. Also, using
the fact $\Gamma(1/2)=\sqrt{\pi}$ and $\Gamma(3/2)=\sqrt{\pi}/2$, we
obtain its global Pad\'{e} approximations [1/2]:
\vskip - 13pt
\begin{equation}\label{eq26}
    E_{\frac{1}{2},\frac{3}{2}}(-x)\approx \frac{\frac{2}{\sqrt{\pi}}+\frac{4-\pi}{\pi-2}x}
    {1+\frac{\sqrt{\pi}}{\pi-2}x+\frac{4-\pi}{\pi-2}x^2}.
\end{equation}

\begin{center}
\includegraphics[width=\columnwidth]{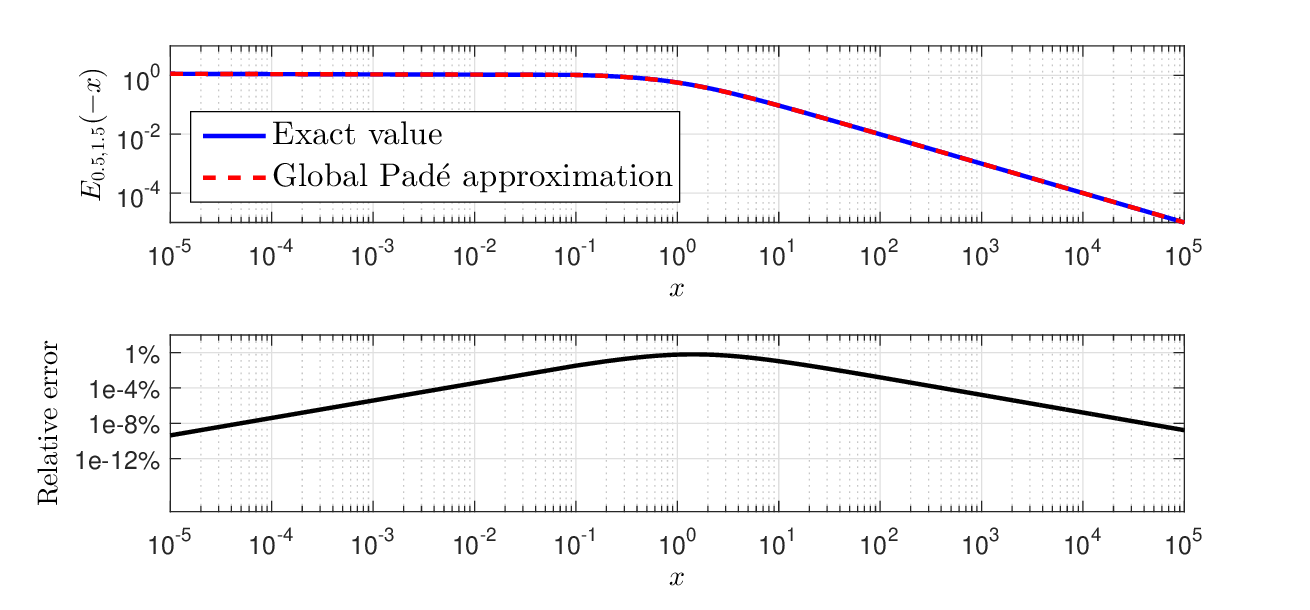}

 \medskip

Fig. 2.1: The top one is the
evolution of the generalized Mittag-Leffler function
$E_{{1}/{2},{3}/{2}}(-x)$ and its global Pad\'{e} approximation [1/2], and the bottom one is the relative error with the peak 0.6389\%.
\end{center}

\begin{remark} 
Numerical simulations in Fig. 2.1 show that (\ref{eq26})
approximates $E_{{1}/{2},{3}/{2}}(-x)$ with relative error
less that $1\%$ for any $x\ge0$.
Also, the approximation (\ref{eq26}) is falling below the exact
value for $x>0$.
\end{remark}

As a special situation, say, $\beta=1$, we reduce (\ref{eq13}) to
\begin{equation}\label{eq14}
   E_{\alpha}(-x)  \approx
   \frac{1+\frac{1}{\Gamma(1-\alpha)q_0^*}x}{1+\frac{q_1^*}{q_0^*}x+\frac{1}{q_0^*}x^2}\,,
\end{equation}
where
\begin{equation}\label{eq25}
\begin{array}{l}
 q_0^*=\dfrac{\frac{\Gamma(1+\alpha)}{\Gamma(1-\alpha)}-
    \frac{\Gamma(1+\alpha)\Gamma(1-\alpha)}
    {\Gamma(1-2\alpha)}}{\Gamma(1+\alpha)\Gamma(1-\alpha)-1},\, ~~
 q_1^*=\dfrac{\Gamma(1+\alpha)-\frac{\Gamma(1-\alpha)}
    {\Gamma(1-2\alpha)}}{\Gamma(1+\alpha)\Gamma(1-\alpha)-1}.
 \end{array}
\end{equation}
It should be pointed out that this special case was studied by
\cite{Atkinson2011}.

\begin{center}
\includegraphics[width=\columnwidth]{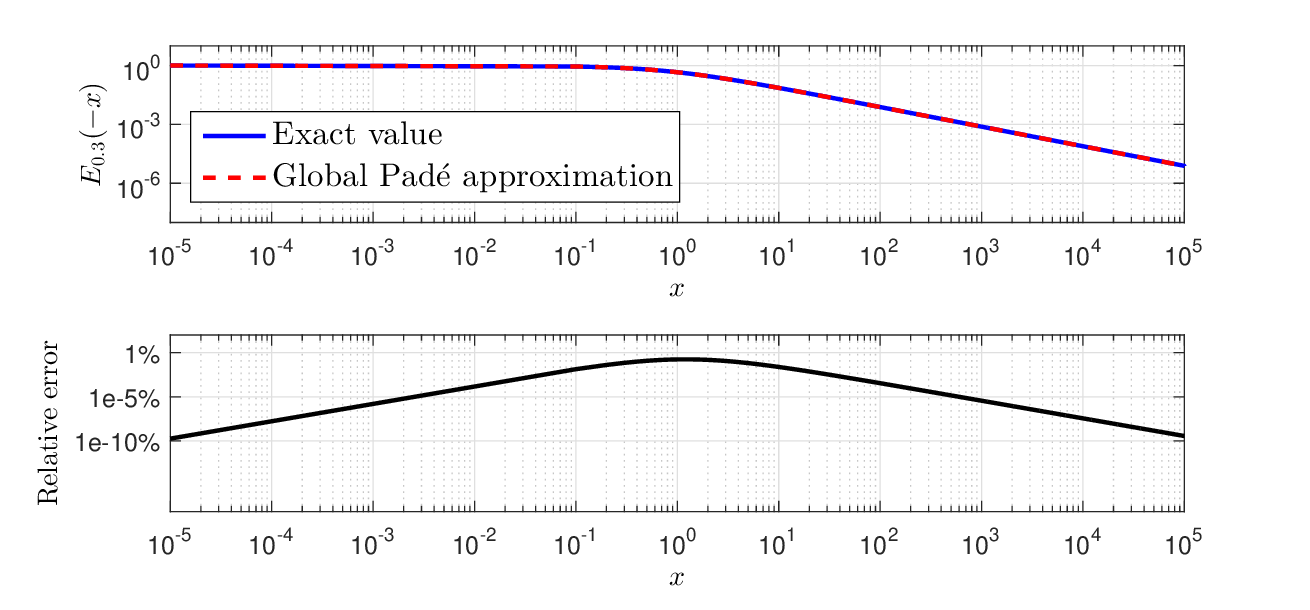}

\medskip

Fig. 2.2: The top one is the
evolution of the generalized Mittag-Leffler function
$E_{0.3}(-x)$ and its global Pad\'{e} approximation [1/2], and the bottom one is the relative error with the peak 0.1744\%.
\end{center}

\begin{remark} 
Numerical simulations in Fig. 2.2 show that (\ref{eq14})
approximates $E_{0.3}(-x)$ with relative error
less that $1\%$ for any $x\ge0$. Also, the approximations
(\ref{eq14}) are smaller than the exact value of $E_{0.3}(-x)$ for
$x\in(0,+\infty)$.
\end{remark} 

When $0<\alpha=\beta<1$, the generalized Mittag-Leffler
function $E_{\alpha,\alpha}(-x)$ admits the two series
\begin{equation}\label{eq15}
\left\{
\begin{aligned}
    \frac{\Gamma(1-\alpha)}{\alpha}x^2E_{\alpha,\alpha}(-x)&=
    \frac{\Gamma(1-\alpha)}{\alpha}x^2
    \sum_{k=0}^{m-3}\frac{(-x)^k}{\Gamma(\alpha+\alpha k)}+\mathcal
    {O}(x^m),\\
        \frac{\Gamma(1-\alpha)}{\alpha}x^2E_{\alpha,\alpha}(-x)&=
    -\frac{\Gamma(1-\alpha)}{\alpha}x^2
    \sum_{k=2}^{n+1}\frac{(-x)^{-k}}{\Gamma(\alpha-\alpha k)}+\mathcal
    {O}(x^{-n}),
\end{aligned}
\right.
\end{equation}
\vskip -4pt \noindent
at $x = 0$ and $x = +\infty$, respectively. The multiplication by
${\Gamma(1-\alpha)}x^2/{\alpha}$ ensures that the first coefficient
of the asymptotic series (\ref{eq15}) is 1.

By the same approach, we can obtain a global Pad\'{e} approximations
[1/2] for $E_{\alpha,\alpha}(-x)$ of the form
\vskip -10pt
\begin{equation}\label{eq17}
    E_{\alpha,\alpha}(-x)\approx \frac{\frac{1}{\Gamma(\alpha)}}
    {1+\frac{2\Gamma(1-\alpha)^2}{\Gamma(1+\alpha)\Gamma(1-2\alpha)}x+
    \frac{\Gamma(1-\alpha)}{\Gamma(1+\alpha)}x^2}\,.
\end{equation}


\begin{center}
\includegraphics[width=\columnwidth]{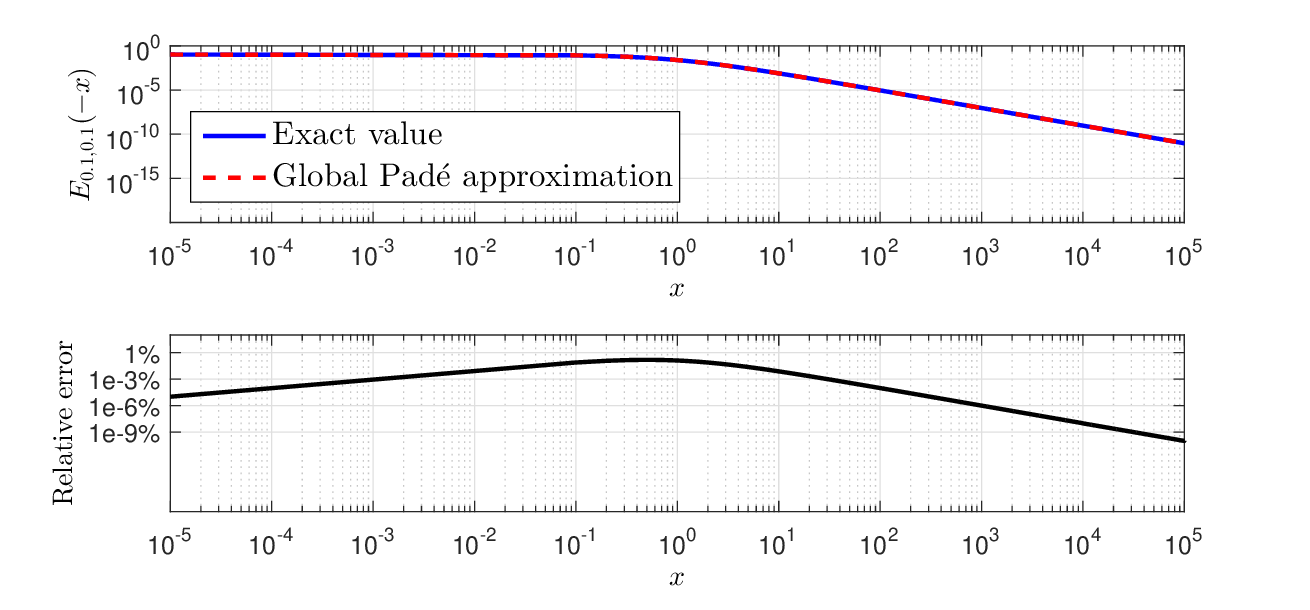}

\medskip

Fig. 2.3: The top one is the
evolution of the generalized Mittag-Leffler function
$E_{0.1, 0.1}(-x)$ and its global Pad\'{e} approximation [1/2], and the bottom one is the relative error with the peak 0.1526\%.
\end{center}

\begin{remark} 
Numerical simulations in Fig. 2.3 show that (\ref{eq17})
approximates $E_{0.1, 0.1}(-x)$ with relative error
less that $1\%$ for any $x\ge0$. Also, the approximations
(\ref{eq17}) are bigger than the exact value of $E_{0.1, 0.1}(-x)$ for
$x\in(0,+\infty)$.
\end{remark}  

\vspace*{-5pt}

Indeed, if $\alpha=1,~\beta>1$, the generalized Mittag-Leffler
function $E_{1,\beta}(-x)$ admits the two series
\vskip -13pt
\begin{equation}\label{eq19}
\left\{
\begin{aligned}
    \Gamma(\beta-1)xE_{1,\beta}(-x)&=
    \Gamma(\beta-1)x
    \sum_{k=0}^{m-2}\frac{(-x)^k}{\Gamma(\beta+k)}+\mathcal
    {O}(x^m),\\
    \Gamma(\beta-1)xE_{1,\beta}(-x)&=
    -\Gamma(\beta-1)x
    \sum_{k=1}^{n}\frac{(-x)^{-k}}{\Gamma(\beta- k)}+\mathcal
    {O}(x^{-n}),
\end{aligned}
\right.
\end{equation}
\vskip -4pt \noindent
at $x = 0$ and $x = +\infty$, respectively. The multiplication by
$\Gamma(\beta-1)x$ ensures that the first coefficient of the
asymptotic series (\ref{eq19}) is 1.

By the same approach, we can obtain a global Pad\'{e} approximations
[1/2] for $E_{1,\beta}(-x)$ of the form
\vskip -13pt
\begin{equation}\label{eq20}
    E_{1,\beta}(-x)\approx \frac{\frac{1}{\Gamma(\beta)}+\frac{1}{\Gamma(\beta+1)}x}
    {1+\frac{2}{\beta}x+
    \frac{1}{\beta(\beta-1)}x^2}\, .
\end{equation}

\begin{center}
\includegraphics[width=\columnwidth]{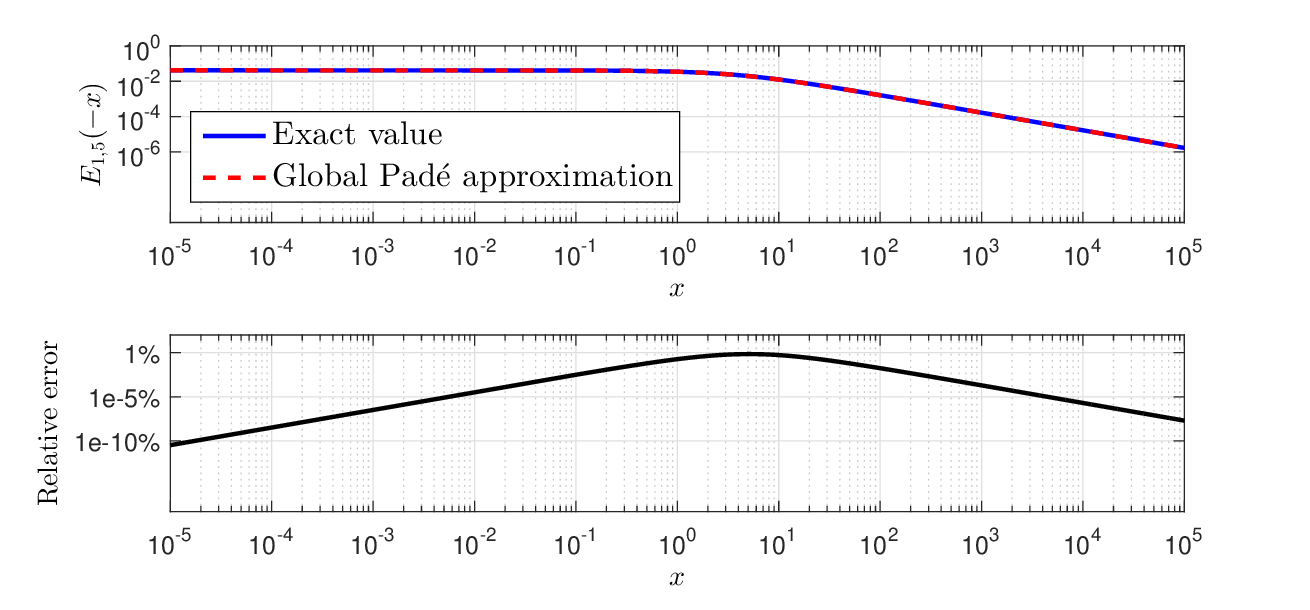}

\medskip

Fig. 2.4: The top one is the
evolution of the generalized Mittag-Leffler function
$E_{1, 5}(-x)$ and its global Pad\'{e} approximation [1/2], and the bottom one is the relative error with the peak
0.6903\%.

\end{center}

\begin{remark} 
Numerical simulations in Fig. 2.4 show that (\ref{eq20})
approximates $E_{1,5}(-x)$ with relative error
less that $1\%$ for any $x\ge0$. Also, the approximation (\ref{eq20}) is falling below the exact
value for $x>0$.
\end{remark}


\begin{remark} 
From (\ref{eq6})-(\ref{eq8}), we can obtain
the approximation with another degree.
In general, higher order approximations, i.e.,
values of $\nu$ greater than 2, can be computed to gain greater accuracy.
However, we find that the proposed approximation [1/2] is better than the Pad\'e approximations
[1/1], [2/3], [3/4] after numerical comparisons.
\end{remark}  

More precisely, we collect the above global Pad\'{e} approximations
[1/2] of the generalized Mittag-Leffler function
$E_{\alpha,\beta}(-x)$, $x\in[0,+\infty)$ for particular cases,
respectively. See Table 1 for more details.

\section{Global Pad\'{e} approximations of $L_{\alpha,\beta}{(x)}$}
\label{sec:3}

\setcounter{section}{3}
\setcounter{equation}{0}\setcounter{theorem}{0}

Hilfer and Seybold \cite{Hilfer2006} introduced and studied the
inverse generalized Mittag-Leffler function $L_{\alpha,\beta}{(x)}$
defined as the solution of the equation \vskip -9pt
\[L_{\alpha,\beta}{(E_{\alpha,\beta}(x))}=x.\]
In the complex plane, the authors have succeeded in determining the
principal sheet of $L_{\alpha,\beta}{(x)}$ provided that three
conditions are satisfied: (1) the function $L_{\alpha,\beta}{(x)}$
is single valued and well defined on its principal sheet; (2) its
principal sheet reduces to the principal sheet of the logarithm
for $\alpha\rightarrow 1$; (3) its principal sheet is a simply
connected subset of the complex plane.

Spencer et al. \cite{Spencer2007} provided also the asymptotic
formulas for the zeros of the function $L_{\alpha,\beta}{(x)}$. In
particular, Atkinson and Osseiran \cite{Atkinson2011} found a global
rational approximation of the inverse Mittag-Leffler function
$-L_{\alpha,1}(x)$ $=-L_{\alpha}(x)$ with $x\in(0,1]$ and $0<\alpha<1$.

It is obvious from the complete monotonicity of the generalized
Mittag-Leffler function $E_{\alpha,\beta}{(-x)}$, described in
(\ref{eq3}) for $x\in[0,+\infty)$, that $E_{\alpha,\beta}{(-x)}$ is
a decreasing and continuous function on the interval $[0,+\infty)$.
This means that the corresponding inverse Mittag-Leffler function
$-L_{\alpha,\beta}(x)$ is well defined on the interval $(0,
1/\Gamma(\beta)]$ using these properties and the fact that
$E_{\alpha,\beta}(0) = 1/\Gamma(\beta)$ by definition.

Our ability to calculate the global Pad\'{e} approximation of
$E_{\alpha,\beta}(-x)$ in previous section
allows us to evaluate also the inverse
generalized Mittag-Leffler function $-L_{\alpha,\beta}(x)$. In fact,
the uniform approximations (\ref{eq13}), (\ref{eq17}), and
(\ref{eq20}) can be inverted to yield a global Pad\'{e}
approximation of the inverse Mittag-Leffler function
$-L_{\alpha,\beta}(x)$ for particular cases $\{0<\alpha<1,
\beta>\alpha\}$, $\{0<\alpha=\beta<1\}$, and $\{\alpha=1,
\beta>1\}$.

For the case $\{0<\alpha<1, \beta>\alpha\}$, by rearranging the
approximation (\ref{eq13}), solving the resulting quadratic
equation, and rearranging, we find that a global Pad\'{e}
approximation of $-L_{\alpha,\beta}(x)$:
\vskip -12pt
\begin{eqnarray}\label{eq21}
    -L_{\alpha,\beta}(x)&\approx&
    \frac{1}{2\Gamma(\beta-\alpha)x}-\frac{q_1}{2} \nonumber\\
    & + & \sqrt{\left(\frac{q_1}{2}-\frac{1}{2\Gamma(\beta-\alpha)x}\right)^2
    -q_0\left(1-\frac{1}{\Gamma(\beta)x}\right)}\ ,
\end{eqnarray}
\vskip -2pt \noindent
where $-L_{\alpha,\beta}(x)|_{x\rightarrow 0^+}=+\infty$ and
$-L_{\alpha,\beta}(x)|_{x\rightarrow 1/\Gamma(\beta)}=0$.

\vskip 3pt

When $\beta=1$, the approximation (\ref{eq21}) reduces to
\vskip -12pt
\begin{eqnarray}\label{eq22}
    -L_{\alpha}(x) &\approx&
    \frac{1}{2\Gamma(1-\alpha)x}-\frac{q_1^*}{2}\nonumber\\
    & + & \sqrt{\left(\frac{q_1^*}{2}-\frac{1}{2\Gamma(1-\alpha)x}\right)^2
    -q_0^*\left(1-\frac{1}{x}\right)}\ ,
\end{eqnarray}
\vskip -2pt \noindent
where $-L_{\alpha}(x)|_{x\rightarrow 0^+}=+\infty$ and
$-L_{\alpha}(x)|_{x\rightarrow 1}=0$.

When $0<\alpha=\beta<1$, we also have a global Pad\'{e}
approximation of $-L_{\alpha,\alpha}(x)$:
\vskip -13pt
\begin{eqnarray}\label{eq23}
    -L_{\alpha,\alpha}(x) &\approx&
    -\frac{\Gamma(1-\alpha)}{\Gamma(1-2\alpha)x}\nonumber\\
    & + &\sqrt{\frac{\Gamma(1-\alpha)^2}{\Gamma(1-2\alpha)^2x^2}
    -\frac{\Gamma(1+\alpha)}{\Gamma(1-\alpha)}\left(1-\frac{1}{\Gamma(\alpha)x}\right)}\ ,
\end{eqnarray}
\vskip -4pt \noindent
where $-L_{\alpha,\alpha}(x)|_{x\rightarrow 0^+}=+\infty$ and
$-L_{\alpha,\alpha}(x)|_{x\rightarrow 1/\Gamma(\alpha)}=0$.

Similarly, if $\alpha=1,\beta>1$, we obtain a global Pad\'{e}
approximation of $-L_{1,\beta}(x)$:
\vskip -14pt
\begin{eqnarray*}\label{eq24}
    -L_{1,\beta}(x)&\approx&
    \frac{1}{2\Gamma(\beta-1)x}-\beta+1\nonumber\\
    &&+\sqrt{\left(\beta-1-\frac{1}{2\Gamma(\beta-1)x}\right)^2
    -\beta(\beta-1)\left(1-\frac{1}{\Gamma(\beta)x}\right)},
\end{eqnarray*}
\vskip - 3pt \noindent
where $-L_{1,\beta}(x)|_{x\rightarrow 0^+}=+\infty$ and
$-L_{1,\beta}(x)|_{x\rightarrow 1/\Gamma(\beta)}=0$.

\smallskip

More precisely, we collect the above global Pad\'{e} approximations
of the inverse generalized Mittag-Leffler
function $-L_{\alpha,\beta}(x)$, $x\in(0,1/\Gamma(\beta)]$ for
particular cases, respectively. See Table 1 for more
details.


\begin{sidewaystable} \label{table}
\bigskip
\vspace*{12cm}
\renewcommand\arraystretch{2}
\begin{center}
 \smallskip
  \begin{threeparttable}
\begin{tabular}{c|c|c}\hline
 Parameters & Function & Global Pad\'e approximation  \\\hline
\multirow{2}[4]{*}{$0<\alpha<1,\beta>\alpha$}&$E_{\alpha,\beta}(-x)$&$
   \frac{\frac{1}{\Gamma(\beta)}+\frac{1}{\Gamma(\beta-\alpha)q_0}x}{1+\frac{q_1}{q_0}x+\frac{1}{q_0}x^2}$\\\cline{2-3}
    &$-L_{\alpha,\beta}(x)$&$
       \frac{1}{2\Gamma(\beta-\alpha)x}-\frac{q_1}{2}
    +\sqrt{\left(\frac{q_1}{2}-\frac{1}{2\Gamma(\beta-\alpha)x}\right)^2
    -q_0\left(1-\frac{1}{\Gamma(\beta)x}\right)}$\\\hline
\multirow{2}[4]{*}{$0<\alpha<1,\beta=1$}&$E_{\alpha}(-x)$&
$\frac{1+\frac{1}{\Gamma(1-\alpha)q_0^*}x}{1+\frac{q_1^*}{q_0^*}x+\frac{1}{q_0^*}x^2}$\\\cline{2-3}
&$-L_{\alpha}(x)$& $\frac{1}{2\Gamma(1-\alpha)x}-\frac{q_1^*}{2}
    +\sqrt{\left(\frac{q_1^*}{2}-\frac{1}{2\Gamma(1-\alpha)x}\right)^2
    -q_0^*\left(1-\frac{1}{x}\right)}$\\\hline
\multirow{2}[4]{*}{$0<\alpha=\beta<1$}& $E_{\alpha,\alpha}(-x)$ &
$\frac{\frac{1}{\Gamma(\alpha)}}
    {1+\frac{2\Gamma(1-\alpha)^2}{\Gamma(1+\alpha)\Gamma(1-2\alpha)}x+
    \frac{\Gamma(1-\alpha)}{\Gamma(1+\alpha)}x^2}$ \\\cline{2-3}
& $-L_{\alpha,\alpha}(x)$ & $
-\frac{\Gamma(1-\alpha)}{\Gamma(1-2\alpha)x}
    +\sqrt{\frac{\Gamma(1-\alpha)^2}{\Gamma(1-2\alpha)^2x^2}
    -\frac{\Gamma(1+\alpha)}{\Gamma(1-\alpha)}\left(1-\frac{1}{\Gamma(\alpha)x}\right)}$
    \\\hline
\multirow{2}[4]{*}{$\alpha=1,\beta>1$} & $E_{1,\beta}(-x)$ &
$\frac{\frac{1}{\Gamma(\beta)}+\frac{1}{\Gamma(\beta+1)}x}
    {1+\frac{2}{\beta}x+
    \frac{1}{\beta(\beta-1)}x^2}$ \\\cline{2-3}
    & $-L_{1,\beta}(x)$ & $
\frac{1}{2\Gamma(\beta-1)x}-\beta+1
    +\sqrt{\left(\beta-1-\frac{1}{2\Gamma(\beta-1)x}\right)^2
    -\beta(\beta-1)\left(1-\frac{1}{\Gamma(\beta)x}\right)}$ \\ \hline
\end{tabular}

 \begin{tablenotes}
      \centering
       \item[] The values of $q_0$, $q_1$, $q_0^*$ and $q_1^*$ are listed in
(\ref{eq12}) and (\ref{eq25}).
     \end{tablenotes}

  \end{threeparttable}
  \bigskip
  \caption{Global Pad\'{e} approximations [1/2] of the
   generalized Mittag-Leffler function and its inverse}
\end{center}
\end{sidewaystable}

\vspace*{-4pt}
\section{Concluding discussions}
\label{sec:4}

\setcounter{section}{4}
\setcounter{equation}{0}\setcounter{theorem}{0}

On the basis of the Taylor series and the asymptotic series of
generalized Mittag-Leffler function, we have constructed a global
Pad\'{e} approximation of the function $E_{\alpha,\beta}(-x)$ with
$x\in[0,+\infty)$ for the particular cases $\{0<\alpha<1,
\beta>\alpha\}$, $\{0<\alpha=\beta<1\}$, and $\{\alpha=1,
\beta>1\}$, respectively. These uniform approximations pick up the
initial exponential-type behavior of the generalized Mittag-Leffler
function as well as its asymptotic power laws for large arguments.
Moreover, these approximations were inverted to yield a global
Pad\'{e} approximation of the inverse generalized Mittag-Leffler
function $-L_{\alpha,\beta}(x)$.

Atkinson and Osseiran \cite{Atkinson2011} reported the rational
solution of time-fractio\-nal diffusion equation under Caputo
definition by using the global Pad\'{e} approximation of the
Mittag-Leffler function $E_{\alpha}(-x)$ with $0<\alpha<1$. However,
this does not suffice to derive the rational solution of the ordinary
and partial time-fractional differential equations in the sense of
Riemann-Liouville. To demonstrate the advantage of our established
results, let us consider the following illustrative examples.

\vskip 3pt

For example, consider the linear fractional differential equation
\vskip -10pt
\begin{equation*}
    {}_{0}D_t^\alpha f(t)+\lambda f(t)=0, \ \ t\ge0, \ \left[{}_{0}D_t^{-\alpha}
f(t)\right]_{t=0}=C_1,
\end{equation*}
\vskip -3pt \noindent
where $0<\alpha<1$, $\lambda>0$, $C_1$ is some constant,
${}_{0}D_t^\alpha$ and ${}_{0}D_t^{-\alpha}$ are the
Riemann-Liouville fractional derivative and integral
of order $\alpha$, respectively. It admits
the following solution \cite{Oldham1974} \vskip -10pt
\begin{equation*}
    f(t)=C_1t^{-\alpha}E_{\alpha,\alpha}\left(-\lambda
t^\alpha\right).
\end{equation*}
Using the approximation (\ref{eq17}), we have the global Pad\'{e}
approximation [1/2]:
\vskip -15pt
\begin{equation*}
    f(t)\approx\frac{C_1}{\Gamma(\alpha)t^\alpha+\frac{2\lambda\Gamma(1-\alpha)^2}{\alpha\Gamma(1-2\alpha)}t^{2\alpha}
+\frac{\lambda^2\Gamma(1-\alpha)}{\alpha}t^{3\alpha}}\,.
\end{equation*}

 We consider another example of the form
\newpage 
\begin{equation*}
    {}_{0}D_t^\alpha g(t)+{}_{0}D_t^\beta g(t)=\delta(t), \ \  t\ge0,
\end{equation*}
where $0<\alpha<\beta<1$ and $\delta(t)$ stands for the the Dirac
delta function. It yields the following solution \cite{Podlubny1999}
\vskip -11pt
\begin{equation*}
    g(t)=\left(C_2+1\right)t^{\beta-1}E_{\beta-\alpha,\beta}\left(-t^{\beta-\alpha}\right),
\end{equation*}
\vskip -3pt \noindent
where $C_2=\left[{}_{0}D_t^\alpha g(t)+{}_{0}D_t^\beta
g(t)\right]_{t=0}$. Using the approximation (\ref{eq13}), we have
the global Pad\'{e} approximation:
\vskip -10pt
\begin{equation*}
    g(t)\approx\frac{\frac{C_2+1}{\Gamma(\beta)}t^{\beta-1}
+\frac{C_2+1}{\Gamma(\alpha)q_0'}t^{2\beta-1-\alpha}}{1+\frac{q_1'}{q_0'}t^{\beta-\alpha}+\frac{1}{q_0'}t^{2(\beta-\alpha)}}\,,
\end{equation*}
\vskip -6pt \noindent
where
\vskip -14pt
$$
 q_0'=\frac{\frac{\Gamma(\beta)^2\Gamma(2\beta-\alpha)}{\Gamma(\alpha)}-
    \frac{\Gamma(\alpha)\Gamma(\beta)\Gamma(2\beta-\alpha)}
    {\Gamma(2\alpha-\beta)}}{\Gamma(\alpha)\Gamma(2\beta-\alpha)-\Gamma(\beta)^2},\ ~~
 q_1'=\frac{\Gamma(\beta)\Gamma(2\beta-\alpha)-\frac{\Gamma(\alpha)\Gamma(\beta)^2}
    {\Gamma(2\alpha-\beta)}}{\Gamma(\alpha)\Gamma(2\beta-\alpha)-\Gamma(\beta)^2}\, .
 $$

From the above discussions, our constructed global Pad\'{e}
approximations of the generalized Mittag-Leffler function are quite
effective in constructing the rational solution to fractional
differential equations across an infinite range of the argument.

\vspace*{-6pt}
\section*{Acknowledgements}

This work was partly supported by the National Natural Science Foundation of China (No. 11271139, 11501216), and the Fundamental Research Funds for the Central Universities (No. 2014ZB0033). The authors thank Prof. Francesco Mainardi for his critical comments and helpful
suggestions. Moreover, we would like to thank the referees and Prof. Richard Martin for carefully reading our manuscript and for giving such
constructive comments which substantially helped improving the quality of the paper.

\vspace*{-3pt}



 \bigskip \smallskip

 \it

 \noindent
$^1$ School of Mathematics, 
South China University of Technology \\
Guangzhou 510640, CHINA \\[4pt]
e-mail: zeng.cb@mail.scut.edu.cn (C. Zeng)
\hfill Received: January 15, 2015 \\[12pt]
$^2$ Mechatronics, Embedded Systems and Automation (MESA) LAB \\
School of Engineering, 
University of California, Merced \\
5200 N Lake Road, Merced, CA 95343, USA \\[4pt]
e-mail: yqchen53@ucmerced.edu (Y.Q. Chen)

\rm

\bigskip  

\vskip 0.05cm  
\hrule width40mm height0.10mm 
\vskip 0.05cm

\bigskip

 \noindent
 Please cite to this paper as published in:

\medskip

\emph{Fract. Calc. Appl. Anal.}, Vol. \textbf{18}, No 6 (2015), pp. 1492--1506,

 DOI: 10.1515/fca-2015-0086

\end{document}